\documentclass[a4paper]{article}

\usepackage{bezier,color}

\usepackage{amsmath,latexsym,amssymb}
\usepackage[german,english]{babel}
\newtheorem{theorem}{Theorem}
\newtheorem{prop}{Proposition}[section]
\newtheorem{proposition}[prop]{Proposition}

\newenvironment{definition}
{\stepcounter{thm}\begin{description}\item[Definition
\thesection.\arabic{thm}:]}
{\end{description}}

\newenvironment{example}
{\stepcounter{thm}\begin{description}\item[Example
\thesection.\arabic{thm}:]}
{\end{description}}

\newenvironment{remark}
{\stepcounter{thm}\begin{description}\item[Remark
\thesection.\arabic{thm}:]}
{\end{description}}

\newenvironment{notation}
{\stepcounter{thm}\begin{description}\item[Notation
\thesection.\arabic{thm}:]}
{\end{description}}

\newenvironment{acknowledgement}
{\stepcounter{thm}\begin{description}\item[Acknowledgement
\thesection.\arabic{thm}:]}
{\end{description}}

\newcommand{\proof}[1][]{{\it Proof#1: }}
\newcommand{\qed}[1][3mm]{\hspace*{\fill} $\Box$ \vspace{#1}}

\newcommand{\s}{\sigma}
\newcommand{\br}{\mathrm{Br}}
\newcommand{\diffo}{\mathrm{Diff}^o}
\newcommand{\AAA}{{\mathbf A}}
\newcommand{\CC}{{\mathbf C}}
\newcommand{\PP}{{\mathbf P}}
\newcommand{\tto}{\longrightarrow}
\newcommand{\del}{\partial}
\newcommand{\bfami}{{\cal B}}
\newcommand{\cfami}{{\cal C}}\newcommand{\dfami}{{\cal D}}
\newcommand{\ffami}{{\cal F}}\newcommand{\mfami}{{\cal M}}
\newcommand{\ofami}{{\cal O}}\newcommand{\pfami}{{\cal P}}
\newcommand{\ufami}{{\cal U}}
\newcommand{\ON}[1]{\mathrm{#1}}
\newcommand{\inv}{^{-1}}
\newcommand{\inj}{\hookrightarrow}
\newcommand{\bfrak}{\mathbf b}
\newcommand{\brS}{\ON{Br}^{\text{\sf s}}}

\newcommand{\eqed}{\hspace*{\fill}\qed}

\begin{document}

\title{Bifurcation braid monodromy of plane curves}

\author{Michael L\"onne}

\maketitle

\abstract{
We consider spaces of plane curves in the setting of algebraic geometry
and of singularity theory.
On one hand there are the complete linear systems, on the other we
consider unfolding spaces of bivariate polynomials of Brieskorn-Pham type.
\\
For suitable open subspaces we can define the
\emph{bifurcation braid monodromy} taking values in the Zariski resp.\
Artin braid group. 
In both cases we give the generators of the image.
\\
These results are compared with the corresponding 
geometric monodromy. It takes values in the mapping class
group of braided surfaces. Our final result gives a precise statement
about the interdependence of the two monodromy maps.
\\
Our study concludes with some implication
with regard to the unfaithfulness of the geometric monodromy (\cite{w})
and the - yet unexploited -
\emph{knotted geometric monodromy}, which takes the ambient
space into account.
}

\section{Introduction}
\label{sec:1}
Algebraic geometry and Singularity theory share their interest in
discriminant complements.
Though they look at different parameter spaces, there is an 
apparent common feature:
A closed subset, the \emph{discriminant}, parametrizes objects which
are special (singular) in some sense,
and thus distinguished from generic objects corresponding to points
in the discriminant complement.

Familiar examples are provided by linear systems
resp.\ versal unfoldings, which contain non-singular and singular elements.
\\

The appropriate tool to study the topology of divisor complements seem
to be braid monodromy maps. They were successfully exploited to
determine fundamental groups. But their domain is naturally
the set of regular values of a projection of the divisor complement.

That is different in the situation considered in \cite{bifbraid} 
where the divisor
complement itself is the domain of an interesting braid monodromy map.
Therefore the natural question was raised, whether spaces of
hypersurfaces in singularity theory and projective geometry
support braid monodromy maps and can be better understood by
their study.
\\
 
To get bifurcation braid monodromy we discard the 
\emph{degeneracy locus} consisting of
those polynomials which are special with respect to a chosen projection
of the domain. 
The remaining \emph{generic} polynomials are 
those on which the induced map is a Morse function. In particular
its \emph{bifurcation set} consists
of only finitely many critical values, their constant cardinality
is determined by the topology.
Of course paths of such points in the target will precisely be the
braids associated to paths of generic polynomials, cf.\ section \ref{bif}.
\\
 
We focus our study to the case of plane projective curves
and their close cousins, unfoldings of singular bivariate polynomials of 
Brieskorn-Pham type $f=x^k+y^{n+1}$.

Our results are expressed in terms of the band generators
$\s_{ij}$ which are natural conjugates of the Artin generators
$\s_{i}=\s_{i\;i+1}$ and which
can be identified with the following braid diagrams:


\begin{figure}
\unitlength=8mm
\begin{picture}(15,4.5)(-1.1,0)


\put(1,4){\makebox(0,-0.051){$n$}}
\put(2,4){\makebox(0,0){$n-1$}}
\put(4,4){\makebox(0,0){$j+1$}}
\put(5,4){\makebox(0,0){$j$}}
\put(6,4){\makebox(0,0){$j-1$}}
\put(9,4){\makebox(0,0){$i+1$}}
\put(10,4){\makebox(0,0){$i$}}
\put(11,4){\makebox(0,0){$i-1$}}
\put(13,4){\makebox(0,0){$2$}}
\put(14,4){\makebox(0,0){$1$}}

\bezier{100}(1,1)(1,2)(1,3)
\bezier{100}(2,1)(2,2)(2,3)
\bezier{100}(4,1)(4,2)(4,3)
\bezier{100}(11,1)(11,2)(11,3)
\bezier{100}(13,1)(13,2)(13,3)
\bezier{100}(14,1)(14,2)(14,3)

\bezier{100}(5,1)(5,1.5)(5.5,1.6)
\bezier{100}(9.5,2.4)(10,2.5)(10,3)
\bezier{100}(5.5,1.6)(6,1.7)(7,1.9)
\bezier{100}(8,2.1)(9,2.3)(9.5,2.4)

\bezier{60}(5,3)(5,2.5)(5.5,2.4)
\bezier{190}(5.5,2.4)(7.5,2)(9.5,1.6)
\bezier{60}(9.5,1.6)(10,1.5)(10,1)

\bezier{50}(6,1)(6,1.2)(6,1.5)
\bezier{50}(6,1.8)(6,2)(6,2.2)
\bezier{50}(6,2.5)(6,2.8)(6,3)

\bezier{50}(9,1)(9,1.2)(9,1.5)
\bezier{50}(9,1.8)(9,2)(9,2.2)
\bezier{50}(9,2.5)(9,2.8)(9,3)

\put(3,2){\makebox(0,0){$\cdots$}}
\put(12,2){\makebox(0,0){$\cdots$}}
\put(7.5,2.5){\makebox(0,0){$\cdots$}}

\end{picture}
\begin{center}
\caption{A generator $\s_{ij}$ of the 'dual' or}
BKL presentation of the
braid group.
\end{center}
\end{figure}

\begin{theorem}
\label{bbmbp}
Suppose the singular polynomial $f=y^{n+1}+x^k$ of Brieskorn-Pham
type is considered with respect to the projection along the $y$-coordinate.
Then in $\br_{nk}$ the conjugation class of the bifurcation braid
monodromy group is represented by
$$
\left\langle \s_{ij}^{m_{ij}}\,\left|\,m_{ij}=\left\{
\begin{array}{ll}
1 & \text{if }i\equiv j\mod n,\\
3 & \text{if }i\equiv j\pm1,\{i,j\}\not\equiv\{0,1\}\mod n,\\
2 & \text{else }
\end{array}
\right.\right.
\right\rangle.
$$
\end{theorem}

\begin{center}
\begin{figure}
\unitlength=1.2mm
\begin{picture}(110,28)(-6.2,10)
\put(0,0){
\begin{picture}(50,40)(-25,-15)
 
\put(20,0){\circle{.8}}
\put(17,10){\circle{.8}}
\put(10,17){\circle{.8}}
\put(0,20){\circle{.8}}
\put(-10,17){\circle{.8}}
\put(-17,10){\circle{.8}}
\bezier{100}(20,0)(18.5,5)(17,10)
\bezier{100}(17,10)(13.5,13.5)(10,17)

\bezier{100}(0,20)(-5,18.5)(-10,17)
\bezier{100}(-10,17)(-13.5,13.5)(-17,10)

\bezier{230}(20,0)(5,8.5)(-10,17)
\bezier{130}(17,10)(8.5,15)(0,20)
\bezier{230}(17,10)(0,10)(-17,10)
\bezier{200}(10,17)(0,17)(-10,17)

\bezier{220}(-10,17)(-13.5,10.5)(-17,4)
\bezier{308}(17,10)(1.2,6)(-14.6,2)

\put(-19.75,2){\circle{.2}}
\put(-20,0){\circle{.2}}
\put(-19.75,-2){\circle{.2}}

\end{picture}
}

\put(50,0){
\begin{picture}(50,40)(-25,-15)
 
\put(20,0){\circle{.8}}
\put(17,10){\circle{.8}}
\put(10,17){\circle{.8}}
\put(0,20){\circle{.8}}
\put(-10,17){\circle{.8}}
\put(-17,10){\circle{.8}}
\bezier{380}(20,0)(3.5,5)(-17,10)
\bezier{200}(20,0)(15,8.5)(10,17)

\bezier{100}(0,20)(5,18.5)(10,17)
\bezier{200}(0,20)(-8.5,15)(-17,10)

\bezier{50}(-18.5,5)(-17.75,7.5)(-17,10)
\bezier{200}(-14,3.4)(-2,10.2)(10,17)

\bezier{30}(20,0)(10,10)(0,20)
\bezier{30}(17,10)(3.5,13.5)(-10,17)
\bezier{30}(10,17)(3.5,13.5)(-17,10)

\bezier{22}(0,20)(-7.5,12.5)(-15,5)
\bezier{40}(-13,0)(3.5,0)(20,0)
 
\put(-19.75,2){\circle{.2}}
\put(-20,0){\circle{.2}}
\put(-19.75,-2){\circle{.2}}

\end{picture}
}
\end{picture}
\begin{center}
\caption{
Both figures depict part of the branch points of $y^4-4y+3x^4$}
at the zeroes
of $x^{12}=1$, so $k=4,l=3$.
The arcs shown on the left \\
hand side correspond
to the twists $\s_{ij}$
with $m_{ij}=3$, those on the right\\
hand side correspond to twists
$\s_{ij}$ with $m_{ij}=2$ and $m_{ij}=1$ (dotted).
\end{center}
\end{figure}
\end{center}

The corresponding result for linear systems of plane curves can be
obtained for the open set of plane curves which are transversal to the line
at infinity and which do not contain the center of projection $(0:1:0)$:

\begin{theorem}
\label{bbmpc} 
The bifurcation braid monodromy group of plane projective curves of
degree $d$ is in the conjugation class of the subgroup
of $\br_{d(d-1)}$ generated by the following elements:
\begin{enumerate}
\item
$\s_{ij}$, if $i\equiv j\mod d-1$,\\[-3mm]
\item
$\s_{ij}^3$, if $i\equiv j\pm1,\{i,j\}\not\equiv\{0,1\}\mod d-1$,\\[-3mm]
\item
$\s_{ij}^2$, if $i,j$ not as above, $i)$ or $ii)$,\\[-3mm]
\item
$\s_1\s_2\cdots\s_{d^2-d-1}\s_{d^2-d-1}\cdots\s_2\s_1$ and
its conjugates by powers of \\
\hspace*{\fill}$\s_{d^2-d-1}\s_{d^2-d-2}\cdots\s_2\s_1$.
\end{enumerate}
\end{theorem}

Both bifurcation braid monodromy groups are in fact isomorphic
to a group of mapping classes, see Prop.\ \ref{isom} for a proof
in case of bivariate polynomials.
These mapping classes are obtained as the natural images of
a \emph{braided geometric monodromy} to be defined in Section
\ref{geo}.
\\

This isomorphism is another instance of the close connection between
algebraic geometry and low dimensional topology,
which is witnessed also by
\begin{enumerate}
\item
the isomorphism induced by geometric monodromy
between the
(orbifold) fundamental group of moduli spaces of curves and
the mapping class group of the corresponding topological surface,
\item
the isomorphism between
the fundamental group of the space of simple polynomials and
the braid group, see Section \ref{bif},
\item
geometric monodromy of plane curves, which induces an injection of the
fundamental group of the discriminant complement of polynomials
of type  $A_n$ and $D_n$
into the mapping class group \cite{pv}.
\end{enumerate}

Our ongoing projects aim at a corresponding result in the absence
of a projection map. Then there is a kind of 
\emph{knotted geometric monodromy}
with range in the mapping classes 
of pairs consisting of an ambient space and the embedded hypersurface.
\\

In the basis case that the family of ambient spaces is trivial
there is another natural candidate for the range of the geometric
monodromy, the fundamental group of higher dimensional
configuration spaces.
Their study was proposed by Dolgachev and Libgober \cite{dl} 
as the topological counterpart of spaces of algebraic submanifolds,
e.g.\ smooth projective plane curves in $\PP^2$.

The appropriate topological space should contain all topological
submanifolds isotopic to a smooth curve $C_d$ of degree $d$.
It can be 
identified as a coset space for the group
$\diffo(\PP^2)$ of diffeomorphisms of $\PP^2$ isotopic to the identity
with respect to the subgroup $\diffo(\PP^2,C_d)$ 
of diffeomorphisms which induce a diffeomorphism
of $C_d$ to itself.

This coset space is the natural topological `configuration space' in higher dimensions
$$
F_{C_d}[\PP^2]\quad=\quad \diffo(\PP^2)/\diffo(\PP^2,C_d)
$$
in analogy to $F_d[S^2]=\diffo(S^2)/\diffo(S^2,\{p_1,...,p_d\})$.

The corresponding quotient map is a fibration which gives rise to a homotopy exact sequence
$$
\pi_1\diffo(\PP^2)\tto\pi_1(F_{C_d}[\PP^2])\tto\pi_0\diffo(\PP^2,C_d)\tto 1
$$
where of course the middle group should be called the `generalised' braid group of
algebraic curves in $\PP^2$.
\\

This raises a lot new questions, about the relation between
the knotted mapping class group and the fundamental group of
higher dimensional configuration spaces, and the respective
geometric monodromy maps.

But with the results of this paper it may be conceivable to get
hold on injectivity and surjectivity properties of these monodromy maps.

\section{Singularity theory}
\label{singtheo}

\newcommand{\prx}{q_x}

Let us first briefly review some basic notions of singularity theory.
We restrict our attention to the case of bivariate polynomials from the beginning.
Note that a rigorous treatment would demand the language of germs,
but for the sake of clarity we will naively speak of
polynomials, (plane) curves and affine spaces.

\begin{definition}
A holomorphic function $f$ defined in a neighbourhood of $0\in \CC^2$
defines a singular curve, if $0\in\CC^2$ is a critical point of $f$ with
critical value $0\in\CC$,
$$
f(0) \quad = \quad \del_x f (0) \quad = \quad \del_y f(0) \quad = \quad 0.
$$

\noindent
Two singular functions are called \emph{equivalent}, if they differ by
a change of coordinates only.
\end{definition}

We are also interested in a more restricted equivalence with respect to
a linear projection 
$$
\prx:\quad \CC^2\to \CC, \quad (x,y) \mapsto x.
$$

\begin{definition}
Two singular functions are called \emph{equivalent rel $\prx$},
if they differ by a holomorphic change $\phi$ of coordinates only, 
which fits into
a commutative diagram with a suitable biholomorphic $\psi$:
\begin{eqnarray*}
\CC^2 & \stackrel{\phi}{\tto} & \CC^2 \\
\prx\downarrow\, &  & \,\downarrow \prx \\
\CC & \stackrel{\psi}{\tto} & \CC
\end{eqnarray*}
\end{definition}

The concept of semi-universal unfolding gets hold of all 
local perturbations
of $f$, at least up to equivalence, resp.\ equivalence rel $\prx $.

Suppose now that $f$ and $f\big|_{x=0}$ are isolated singularities.
In that case,
the semi-universal unfolding rel.\ $\prx $ associated to $f$
is given by a function $F$.
It is determined by the respective equivalence class of $f$ up 
to non-canonical isomorphism.
$$
F: \CC^2\times \CC^{\mu+\mu'} \tto \CC,
$$
where $\mu$ is the Milnor number of $f$ 
and $\mu'$ the Milnor number of $f\big|_{x=0}$.

The following \emph{bifurcation diagram} displays the essential
objects and maps for our set-up:
$$
\begin{array}{cclcccl}
x,y,u && \CC^{\mu+\mu'+2} & \supset & \cfami & \\[2mm]
\downarrow && \downarrow \\[2mm]
x,u && \CC^{\mu+\mu'+1} & \supset & \bfami & \\[2mm]
\downarrow && \downarrow \\[2mm]
u && \CC^{\mu+\mu'} & \supset & \dfami
\end{array}
$$

In this diagram we placed some emphasis on the family of plane curves
$\cfami$, the zero set of $F$,
on the \emph{branch divisor}
\begin{eqnarray*}
\bfami & := & \{(x,u)\,|\,F_{x,u} :\, y\mapsto F(x,y,u)
\text{ has singular zero levelset}\} \\
& \phantom{:}= & \{(x,u)\,|\,F_{x,u} :\, y\mapsto F(x,y,u)
\text{ has multiple roots}\} 
\end{eqnarray*}
and on the \emph{degeneracy locus} 
$$
\dfami:=\{u\,|\,\cfami_u \text{ is singular or } 
\prx|_{\cfami_u} \text{ is not a Morse function}\}
$$

\noindent
We note the following features:
\begin{itemize}
\item
$\cfami\to\CC^{\mu+\mu'+1}$ is a finite map with branch locus $\bfami$,
\item
$\bfami\to\CC^{\mu+\mu'}$ is a finite map with branch locus $\dfami$,
\item
$\bfami$ is the zero set of a monic polynomial $p$ 
of degree ${\mu+\mu'}$ in $x$\\ with
coefficients in $\CC[u]$.
\item
$\dfami$ is the locus of parameters such that the corresponding
monic polynomial $p_u$ has a multiple root.
\end{itemize}

In particular there is a well-defined Lyashko-Looijenga map on the
complement of the degenracy locus
$$
U_{(f)} \quad :=\quad
\CC^{\mu+\mu'} - \dfami \quad \tto \quad \CC[x],
\quad
u \quad \mapsto \quad p_u,
$$
\label{unotation}
which maps to monic univariate polynomials 
of degree $\mu+\mu'$ with simple roots only.

\section{Braid monodromy maps and groups}
\label{bif}

A braid monodromy in general is the map on fundamental groups
induced from a topological map on a suitable space to a space
which has a braid group as its fundamental group.
Here we are only interested in the braid group of the plane and the sphere,
i.e.\ the fundamental groups of the associated configuration spaces.

The configuration space $U_d$ of $d$ points in $\CC$ is naturally
an open algebraic subset of the affine space $\AAA_d$ of monic
univariate polynomials of degree $d$. 
Polynomials in $U_d$ are characterized by
the property that they have simple roots only.

\begin{proposition}[\cite{ar}]
The fundamental group of the open subset
$U_d$ is isomorphic to the (planar) \emph{braid group}
$\br_d$.
It is finitely presented by generators  $\s_i$, $1\leq i < d$
and by relations
\begin{enumerate}
\item
$\s_i\s_j \,= \, \s_j\s_i,\quad$ if $|i-j|>1$, $1\leq i,j < d$,\\[-3mm]
\item
$\s_i\s_{i+1}\s_i \,=\, \s_{i+1}\s_i\s_{i+1},\quad$ if $1\leq i < d-1$,
\end{enumerate}
\end{proposition}

The configuration space of $d$ points on $\PP^1$, 
the topological two-sphere,
is naturally an open algebraic subset of the projective space associated
to the vector space  $V_d=\ON{Sym}^d\CC^2$ 
of homogenous polynomials of degree $d$ in two variables.

\begin{proposition}[\cite{z}]
The fundamental group of the open set in 
$\mathbb P H^0(\ofami_{\mathbb P^1}(d))\cong \mathbb P V_d$,
which consists of homogeneous polynomials with simple roots only,
is isomorphic to the \emph{spherical braid group}
$\ON{Br}^{\text{\sf s}}_d$. 
It is finitely presented by generators $\s_i$, $1\leq i < d$
and by relations
\begin{enumerate}
\item
$\s_i\s_j \,= \, \s_j\s_i,\quad$ if $|i-j|>1$, $1\leq i,j < d$,\\[-3mm]
\item
$\s_i\s_{i+1}\s_i \,=\, \s_{i+1}\s_i\s_{i+1},\quad$ if $1\leq i < d-1$,\\[-3mm]
\item
$\s_1^{}\cdots\s_{d-2}^{}\s_{d-1}^2\s_{d-2}^{}\cdots\s_{1}^{} \, = \, 1$.
\end{enumerate}
\end{proposition}

Given now any map to $\AAA_d$ the
restriction to the preimage of $U_d$ induces a map from the 
fundamental group of the preimage to $\br_d$, the fundamental group
of the image.
(There is a close analog in the spherical case, of course.)

In our set-up we may just look at a family of bivariate polynomials
 $f_u(x,y)$ which maps to a family of monic univariate polynomials
 $p_u(x)$ obtained as the resultant with respect to $y$ of $f$ and
 its derivative $\del_y f$ with respect to $y$, 
 $$
 p_u(x) \quad = \quad  \ON{res_y}(f, \del_y f).
 $$
 The induced map on fundamental groups is called the 
 \emph{bifurcation braid monodromy}.
\\

In particular, we can apply these considerations to the 
Lyashko-Looijenga map of the last section. 
Note again the point stressed there, that the complement $U_{k;n}$
of the degeneracy locus $\dfami$ is mapped to the discriminant complement
$U_{nk}$.

\begin{definition}
The \emph{bifurcation braid monodromy group} of a bivariate polynomial
$f$ is the image of the bifurcation braid monodromy for the versal
unfolding of $f$ relative $\prx $.
\end{definition}

\begin{example}
\label{fullbraid}
The bifurcation braid monodromy of any generic polynomial deformation
of a function equivalent to
$y^2-x^k$ rel.\ $\prx $ is the full braid group $\br_k$.

To see this it suffices to note that a versal family is parametrized
by $\AAA_k$. More precisely the complement of $\dfami$ is $U_k$
and the Lyashko-Looijenga map the identity in this case.
\end{example}

\begin{example}
\label{Anbraids}
The bifurcation braid monodromy of the function
$f=y^{n+1} + x$ rel.\ $\prx $ is in the conjugacy class of the subgroup
of $\br_n$ generated by
$$
\s_i^3,\s_{i,j}^2,|i-j|\geq2.
$$

In fact the bifurcation diagram for this function
is a smooth pull-back of the discriminant diagram of the function $t^{n+1}$.
In terms of an unfolding $F$ and a truncated
unfolding $F'$ of $t^{n+1}$
$$
F(t, u_{n-1},\dots, u_1, u_0) \:=\: F'(t, u_{n-1},\dots, u_1) + u_0
 \:=\:  t^{n+1} + \sum_{i=0}^{n-1} u_i t^i.
$$
the discriminant diagram is given by
$$
\begin{array}{cclcccl}
u=(u_{n-1},\dots, u_1, u_0) && \CC^n & \supset & \dfami & 
= & \{u\:|\:F_u \text{ has singular zero-level }\}\\[2mm]
\downarrow && \downarrow \\[2mm]
u'=(u_{n-1},\dots, u_1) && \CC^{n-1} & \supset & \bfami &
= & \{u'\,|\,F'_{u'} \text{ is not a Morse function}\}.
\end{array}
$$
Then we can refer to the corresponding claim for the braid monodromy
group of the discriminant diagram, which 
is shown in \cite[section 4]{phamcurve} and relies on \cite{loo}
and \cite{cw}.
\end{example}

\section{Computation of bifurcation braid monodromy}
\label{proof}

In this section we want to give the outline of a proof of 
Theorem \ref{bbmbp}.
The following proposition has been proved with more detail in \cite{hurwitz}
in the special case of $l=2$.

\begin{proposition}
\label{hurwitz}
The bifurcation braid monodromy group of a plane curve germ
$y^{n+1}+x^k$ projected by
the $x$-coordinate is the subgroup of $\br_{nk}$ generated by
\begin{enumerate}
\item
$\s_{ij}$, if $i\equiv j\mod n$,\\[-3mm]
\item
$\s_{i,i+1}^3$, if there is $s\equiv 0 \mod n$ such that $s< i\leq n+s-1$,\\[-3mm]
\item
$(\s_{ij}')^2$, if there is $s\equiv 0 \mod n$ such that $s< i,j\leq n+s, |i-j|>1$,
\end{enumerate}
where
$\s_{ij}'=\s_{i,i+1}\inv\cdots\s_{j-2,j-1}\inv\s_{j-1,j}\s_{j-2,j-1}\cdots\s_{i,i+1}$.
\pagebreak[4]
\end{proposition}
 
\unitlength=1.2mm

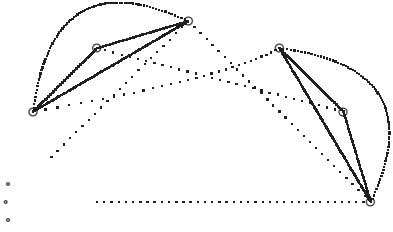
\begin{center}
\begin{picture}(50,30)(-25,-5)
 
\put(20,0){\circle{.8}}
\put(17,10){\circle{.8}}
\put(10,17){\circle{.8}}
\put(0,20){\circle{.8}}
\put(-10,17){\circle{.8}}
\put(-17,10){\circle{.8}}
\bezier{100}(20,0)(18.5,5)(17,10)
\bezier{100}(17,10)(13.5,13.5)(10,17)
\bezier{200}(20,0)(15,8.5)(10,17)

\bezier{100}(0,20)(-5,18.5)(-10,17)
\bezier{100}(-10,17)(-13.5,13.5)(-17,10)
\bezier{200}(0,20)(-8.5,15)(-17,10)

\bezier{80}(20,0)(27,14.9)(10,17)
\bezier{80}(0,20)(-14.9,27)(-17,10)

\bezier{30}(20,0)(10,10)(0,20)
\bezier{30}(17,10)(3.5,13.5)(-10,17)
\bezier{30}(10,17)(3.5,13.5)(-17,10)

\bezier{22}(0,20)(-7.5,12.5)(-15,5)
\bezier{38}(-10,0)(5,0)(20,0)
 
\put(-19.75,2){\circle{.2}}
\put(-20,0){\circle{.2}}
\put(-19.75,-2){\circle{.2}}

\end{picture}
\begin{quote}
\begin{center}
The figure shows part of the branch locus $x^{12}=1$ of
$y^4-4y+3x^4$, \\[0.9mm]
so $k=4,n=3$.
The straight arcs correspond to twists $\s_{ij}$ with\\
$m_{ij}=3$ and $m_{ij}=1$ (dotted),the curved arcs to twists ${\s_{ij}'}^2$ of
the bifurcation braid monodromy group.
\\[5mm]
\end{center}
\end{quote}
\end{center}

\proof
To smooth the argument we note one computational detail in advance.
The branch locus of a polynomial
$f=y^{n+1}-(n+1)p(x)y+nq(x)$ with respect to $\prx $ is obtained
by an elementary elimination:
\begin{eqnarray*}
\del_y f = 0 & \implies &
y^n = p(x) \\
& \stackrel{f=0}{\implies} & 
y p(x) = q(x) \\
& \implies &
y^n p^n(x) = q^n(x) \\
& \stackrel{\del_y f=0}{\implies} &
p^{n+1}(x) = q^n (x) \hspace{2cm}{(*)}
\end{eqnarray*}

In the first step we consider the family 
$f=y^{n+1} - (n+1)u y + n (x^k + v )$ parametrized by
$u,v$.
By definition the bifurcation braid monodromy is induced by the map
$$
u,v \quad \mapsto \quad p_{u,v}(x) = \ON{discr}_y(f)
:=\ON{res}_y(f,\del_y f)
\stackrel{(*)}{=} (x^k + v)^n - u^{n+1}.
$$
(Let us remark that we feel free to rescale the discriminant without further
notice.)

To find the intersection of the $uv$-parameter plane with the open
set of admissible polynomials we have to find the $u,v$ such that
$p_{u,v}$ has a multiple root. Again the elimination of $x$ from
$p_{u,v}$ and $\del_x p$ is quite elementary:
\begin{eqnarray*}
\del_x p = 0 & \implies & x^{k-1}(x^k + v)^{n-1} = 0 \\
& \implies & x^k(x^k + v)^{n-1} = (v - v) (x^k + v)^{n-1}\\
& \stackrel{p=0}{\implies} &
u^{n+1} = v(x^k + v)^{n-1}\\
& \stackrel{p=0}{\implies} &
u^{n(n+1)} = v^n u^{(n-1)(n+1)}\\
& \implies &
u^{(n-1)(n+1)}\big(u^{n+1} - v^n\big) = 0 \hspace{2cm}{(**)}
\end{eqnarray*}

To determine now the braid monodromy we fix a base point
at $(u,v) = (1,0)$. The corresponding branch set is given by
$p_{(1,0)} = x^{nk}-1 = 0$, see $(*)$.

For further use we
number these branch points 
according to increasing $arg$ ending with $\xi_{nk}=1$.
\\

The bifurcation braid monodromy is now defined on the fundamental group
of the complement to $(**)$.
Natural generators are given by a geometric basis on the
line $(1,v)$ punctured where $v^n=1$ and a simple closed path
around $u=0$, see below for details.
\\

On the line $(1,v)$ it suffices to consider paths where $v$ moves along
radial rays from $0$ to a unit root $\xi_{ik}$, $1\leq i \leq n$.
Again from $(*)$ we have for $\lambda\in [0,1]$
$$
(x^k+\lambda \xi_{ik})^n \: = \: 1 \quad \Leftrightarrow \quad
 x^k \: = \: \xi_{jk}-\lambda \xi_{ik} \quad \text{for some } j.
$$
Accordingly the $k$ 
branch points with indices congruent to $i \mod n$ 
converge along radial rays to $0$ while the remaining branch points
stay away from these rays.
\\

The local monodromy at the degeneration points $(1,\xi_{ik})$
can be obtained from the family $g_t = y^2  + x^k - t$.
Its bifurcation divisor $x^k - t  = 0$ is the local model of the
singular branch of the bifurcation divisor of $f$ over $(1,\xi_{ik})$.

Accordingly the monodromy of $g_t$ is mapped to the local
monodromy of $f$ by a transfer map, which identifies a disc
containing the solutions of $x^k = t$ with the disc
containing the $k$ branch points converging to the origin.
\\

We now get to the final path from $(1,v)$ around the line $u=0$.
Let $\rho$ be the solution of $\rho^2=\xi_k$ with positive imaginary part, 
then we can consider the degeneration along
$ (u ,v ) = \big( (1-\lambda)^{\frac n{n+1}}, \lambda \rho \big) $ with $\lambda\in[0,1]$
and bifurcation according to $(*)$:
\begin{eqnarray*}
&&
(x^k + \lambda \rho )^n=(1-\lambda)^{n}\\
& \implies\quad & 
 x^k = (1 - \lambda ) \xi_{jk} - \lambda \rho \quad \text{for some } 
j\in\{1,\dots,n\}.
\end{eqnarray*}
Let $x_{i'}$ denote the solution which moves to $\xi_{i'}$ for $\lambda\to 0$.
Then one may check for $\lambda\to1$
that the argument $arg(x_{i'})$ is strictly increasing
(resp.\ decreasing or constant) 
for $i' \equiv i \mod n$, $0< i < (n+1)/2$ (resp.\ $(n+1)/2<i\leq n$ or $i=(n+1)/2$).

Hence the family degenerates at
$\lambda=1$ only and all branch points are on distinct rays for $\lambda\in[0,1[$.
Moreover we observe that for $\lambda\to 1$ the following $n$-tuples of
branch points merge at $k$ distinct points,
$$
T_s\::=\:\{x_s,x_{s+1},...,x_{s+n-1}\} \quad\text{with}\quad{s\equiv 0 \mod n}.
$$


The local monodromy at the degeneration point $(0,\rho)$
can be obtained from the family $f_t = y^{n+1} -(n+1) t y + n x$.
Its bifurcation divisor $x^n = t^{n+1}$ is the local model of each
singular branch of the bifurcation divisor of $f$ over $(0,\rho)$.

Accordingly each local monodromy of $f$ is obtained by a
transfer map from the monodromy of $f_t$.
The transfer map identifies a disc containing the solutions
of $x^n = t^{n+1}$ with a disc containing an $n$-tupel $T_s$
of branch points converging to a singularity over $(0,\rho)$.
\\

From the monodromy of the special family we may get the monodromy
of the versal family by the principles of versal braid monodromy 
\cite{versal}.
They tell us how to replace the generating braids associated to the
special family by groups of braids, which then generate the
full braid monodromy group of the versal family.
\\

In fact one has to find first the braid monodromy groups of the local
models. Then the transfer maps mentioned above map these group to
subgroups of $\br_{nk}$ which generate the bifurcation braid monodromy
group.

The local model for the degenerations in the line $(1,v)$ is given
by example \ref{fullbraid}, so for each point $(1,\xi_{ik})$ we have
to transfer a full braid group $\br_{k}$ into $\br_{nk}$.
This is done by a topological disc which contains the $\xi_j$ with
$j \equiv i \mod k$. 
Thus the contribution to the bifurcation braid monodromy group
is given by the half-twists $\s_{ij}$ with $i \equiv j \mod k$.

The local model for the degenerations at $(0,\rho)$ is given by example
\ref{Anbraids}, which models the degeneration of each of the $k$
tuples $T_i$ of $n$ branch points. 
Here we have to transfer the corresponding
braids from example 
\ref{Anbraids} to topological discs around the $T_i$.
In this way we get the remaining braids of the claim.
\eqed

\begin{remark}
The last argument admittedly is incomplete, since the transfer map
in the last case is sensitive to the identification of the disc.
Since this ambiguity does not matter in the proof of our theorem,
we do not stress the point here.
\end{remark}

\proof[ of Theorem \ref{bbmbp}]
We have to compare the two subgroups generated by the set of braids
given in the assertion of Prop.\ \ref{hurwitz}, respectively the set of braids
given in the theorem.

Let us focus first on generating braids supported on a topological disc
around an $n$-tupel $T_s$.
In the proposition these elements are listed in $2)$ and $3)$ with $n'=s$.
As we remarked above, we did not actually prove that these elements
generate the monodromy group.
Rather, we showed that under some identification of discs
generators coincide with the generators of example \ref{Anbraids}.

In the claim of the theorem the elements concerned are $\s_{ij}^2, \s_{ij}^3$ listed in
$2)$ and $3)$ with $s< i< j \leq s+n$.
They are readily seen to coincide with the generators of example \ref{Anbraids}
under a suitable map.

Hence by composition we get a local homeomorphism under
which the elements from the proposition are identified with the
elements from the theorem.

Since discs around the $T_s$ may be chosen disjoint we deduce
the existence of a conjugating braid $\beta$ which induces
the above indentifications simultaneously.
\\

In the second step we consider the generators listed under $1)$ of the
proposition and the theorem.
They are the same and generate the same subgroup $H$ of $\br_{nk}$,
but helas we have now to take the conjugation by
$\beta$ into consideration.

We note that $\beta$ is symmetric in the sense
that it commutes with any rigid rotations, which permutes the $n$-tuples
$T_s$. 
Therefore $\beta$ belongs to the subgroup generated by
$$
\delta_1=\s_1\s_{n+1}\dots\s_{(k-1)n+1}, \dots , 
\delta_{n-1} = \s_{n-1}\s_{2n-1}\dots\s_{nk-1}.
$$
The crucial observation is, that $H$ is invariant under conjugation by
the $\delta$:
In fact $\delta_{i'}$ acts on $\s_{i,j}$ with $i\equiv j \mod n$ as
$$
\delta_{i'} \s_{i,j} \delta_{i'}\inv \quad =Ê\quad
\left\{
\begin{array}{cl}
\s_{i,j}\s_{i+1,j+1} \s_{i,j}\inv & \text{ if } i'\equiv i \mod n \\
\s_{i-1,j-1} & \text{ if } i'\equiv i-1 \mod n \\
\s_{i,j} & \text{ else}
\end{array}
\right.
$$
So at this stage we have proved that the monodromy group
of the proposition conjugated by $\beta$ is contained in the
group generated by the elements listed in the theorem.

To finish, it suffices to show that the all elements of the theorem
not considered till now are in fact redundant.
But this can be shown inductively using for
$i<j$ that
$\s_{i,j+n} \s_{j,j+n} \:= \: \s_{j,j+n} \s_{i,j}.
$ 
{}\hspace*{\fill}\qed

\begin{remark}
The bifurcation braid monodromy is precisely the subgroup of $\br_{nk}$
generated by those powers of the band generators $\s_{ij}$
which stabilise the periodic sequence of transpositions
$$
(12), (23), \cdots, (k\, k+1), \:
(12), (23), \cdots, (k\, k+1), \cdots ,(12), \cdots, (k\, k+1)
$$
of length $nk$ under the Hurwitz action.
This sequence encodes of course the finite branched covering of $\CC$
by the curve $C:=\{f_{(1,0)}=0\}$ via $\prx $.
\end{remark}

\section{Monodromy for spaces of plane projective curves}
\label{proj}

The space of plane projective curves of degree $d$ is given by
$\mathbb P H^0_d=\mathbb P H^0(\ofami_{\mathbb P^2}(d))$.
In analogy to the situation in singularity theory we consider open
subsets of curves which have a generic branching property.
Of course they are open subsets in the discriminant complement
$\ufami_d$ corresponding to the set of smooth curves.

\begin{notation}
For a given point $p_0\in\mathbb P^2$ we have the subset of
smooth curves disjoint to $P_0$ which are generic with respect
to the projection $q:\mathbb P^2-\{P_0\}\to \mathbb P^1$ from $P_0$:
\begin{eqnarray*}
\ffami_d & = & \{ C \in \ufami_d \,|\, P_0\notin C, q|_{C} \text { is Morse } \}
\end{eqnarray*}

An open subset is obtained imposing the condition that
a line $L_0$ containing $P_0$ (say at infinity)
has $d$ simple points of intersection with $C$.
\begin{eqnarray*}
\ffami_d' & = & \{ C\in\ufami_d \,|\, C\in\ffami_d, \# C\cap L_0 = d \}
\end{eqnarray*}
\end{notation}

\newcommand{\nd}{\mu_d}
\newcommand{\faser}{F'_d}

\begin{remark}
Thanks to the homogeneity of $\PP^2$ these spaces
do not depend on the choice of a projection center and/or line at
infinity.
\\

If we introduce homogeneous coordinates $(x:y:z)$
such that $P_0=(0:1:0)$ and
$L_0 = \{z=0\}$ we can make the following identifications
(with $\nd=d(d+3)/2$):
\begin{eqnarray*}
\AAA^{\nd }  & \quad = \quad
\{  f \in \CC[x,y,z]_d \,\big|\, f(0,1,0)= 1\,\}
\\
\ffami_d & \quad = \quad
\{ f\in \AAA^{\nd} \,\big|\, \ON{discr}_y(f)\, 
\text{Êhas simple roots only} \}
\\
\ffami_d' & \quad = \quad
\{ f\in \AAA^{\nd} \,\big|\, \ON{discr}_y(f)(x,1) \in U_{d(d-1)} \}
\end{eqnarray*}
The complement of $\ffami_d$ in $\AAA^{\nd }$ is the 
\emph{weak degeneracy locus}
$\dfami$, that of $\ffami_d'$ the \emph{degeneracy locus}
$\dfami'$.
\end{remark}

There is pull-back diagram
\begin{eqnarray*}
\faser\quad & \inj & \ffami'_d \:\ni\: f(x,y,z)\\
\downarrow\quad & & \:\downarrow \qquad\quad \downarrow\\
\{y^d+1\}& \inj & U_d \:\ni\: f(1,y,0)
\end{eqnarray*}
which defines $\faser$ as the fibre of the map on the right
hand side over the element $y^d+1$. 

Of course $\faser$
consists of the polynomials, which can be written as a sum of $y^d+x^d$
with at polynomial of degree $d$ that has $z$ as a factor.


\begin{proposition}
\label{sequence}
There is an natural exact sequence of groups
$$
\pi_1(U_{d;d}) \tto \pi_1(\ffami'_d) \tto \br_d \to 1,
$$
where $U_{d;d}$ is $U_{(f)}$ on page \pageref{unotation} 
of section \ref{singtheo} with $f=y^d+x^d$.
\end{proposition}

\proof
First we want to apply the Zariski theorem on fundamental groups
of divisor complements as proved by Bessis \cite[section 2]{b}.
Consider the map
\begin{eqnarray*}
\AAA^{\nd } & \: \tto\: & \AAA^{\nd -1}
\\
f &  \mapsto &
f(x,y,1)- f(0,0,1)
\end{eqnarray*}
which is the projection along the coefficient of $z^d$.
With a generic choice of parameters $\alpha_1,\alpha_2,\alpha_3$
the fibre  $F$ over 
$$
\{y^d+(\alpha_1x+\alpha_2)y+x^d+\alpha_3 x\}\quad\in\quad\AAA^{\nd -1}
$$ 
intersects the weak degeneracy locus $\dfami$ transversally.
The second contribution to the degeneracy locus is
the pull-back $\pfami$ 
of the discriminant in $V_d$,
the complement of $U_d$.

Hence the Zariski theorem asserts that the following sequence is exact:
$$
\pi_1(F- F\cap\dfami) \tto \pi_1(\ffami_d') \tto \pi_1(\AAA^{\nd -1}-\pfami') \tto 1.
$$
where $\pfami'$ is the pullback of the discriminant in $V_d$ to
$\AAA^{\nd -1}$.

Since $F- F\cap\dfami$ is a subset of $\faser$, the first map
factors through $\pi_1(\faser)$. With 
$\pi_1(\AAA^{\nd -1}-\pfami')=\pi_1(U_d)= \br_d$
we then get the exact sequence
$$
\pi_1(\faser) \tto \pi_1(\ffami_d') \tto \br_d \tto 1.
$$
It remains to identify
the fundamental group on the left with $\pi_1(U_{d;d})$.
To do so
it is natural to consider $\faser$ as a subspace of 
a trivial unfolding of $U_{d;d}$.
We claim that the induced map on fundamental groups 
is then an isomorphism.
The proof for that claim has to be copied from the arguments in
\cite[section 4]{duke}.
\eqed
\\

In fact with some more care, it could be proved that the sequence is even
short exact.

\begin{definition}
The \emph{bifurcation braid monodromy} of smooth plane curves of
degree $d$ is given by the map
\begin{eqnarray*}
\pi_1(\ffami'_d) & \quad & \tto \quad \br_{d(d-1)}\\
\big( \text{resp. } \: \pi_1(\ffami_d) 
&\quad & \tto \quad \brS_{d(d-1)}
\: \text{in the spherical case.}Ê\big)
\end{eqnarray*}
which is - in both cases - induced by the map on polynomials
\begin{eqnarray*}
f=f(x,y,z) &\quad \mapsto \quad  \ON{res}_y( f , {\partial _y} f ).
\end{eqnarray*}

Of course the spherical braid monodromy group can be
given with the same generators, since $\pi_1(\ffami')\to\pi_1(\ffami)$
is surjective.
\end{definition}

\proof[ of Theorem \ref{bbmpc}]
We can now outline the proof of the second theorem.
\\
By the exact sequence of Prop.\ \ref{sequence}, 
the braid monodromy group 
of plane curves of degree $d$ contains
all braids which belong to the braid monodromy group of the
bivariate Fermat polynomial with exponent $d$.
To get sufficiently many additional braids, we have to find
lifts of generators of $\pi_1(U_d)$ and add their images in $\br_{d(d-1)}$.

Again we investigate a special family of curves (here we work with the
homogeneous form)
$$
f_u \quad = \quad y^d - d(u x^{d-1} +z^{d-1}) y + (d-1) x^d.
$$
By $(*)$ the branch locus for the family is given by
$$
p_u(x,z) \quad = \quad (u x^{d-1} +z^{d-1} )^d - x^{d(d-1)} 
$$
and we can compute again the degeneracy locus 
\begin{eqnarray*}
\del_z p = 0  & \implies &
z^{d-2}( u x^{d-1} + z^{d-1} )^{d-1} = 0 \\
& \implies &
z^{d-1}( u x^{d-1} + z^{d-1} )^{d-1} = (u x^{d-1} - u x^{d-1}) 
( u x^{d-1} + z^{d-1} )^{d-1} \\
& \implies &
( u x^{d-1} + z^{d-1} )^{d} = u x^{d-1} ( u x^{d-1} + z^{d-1} )^{d-1} \\
& \stackrel{p=0}\implies &
x^{d(d-1)} = u x^{d-1} ( u x^{d-1} + z^{d-1} )^{d-1} \\
& \implies &
x^{d^2(d-1)} = u^d x^{d(d-1)} ( u x^{d-1} + z^{d-1} )^{d(d-1)} \\
& \stackrel{p=0}\implies &
x^{d^2(d-1)} = u^d x^{d(d-1)} x^{d(d-1)^2} \\
& \implies &
x^{d^2(d-1)}(u^d -1) = 0 
\end{eqnarray*}
This implies $u^d=1$ since $x=0$ is only a solution together with $z=0$,
which has no geometric meaning.
\\

We consider now the image of our family under the restriction map:
$$
f_u \qquad \mapsto \qquad f_u(1,y,0) \quad = \quad y^d - d u  y + (d-1).
$$
The induced mapping on parameter spaces $\CC \to \AAA_d$
is transversal to the discriminant, hence surjective on
fundamental groups.
Thus it suffices to find braids associated to a geometric basis of 
paths in the $u$-parameter plane punctured at $u^d=1$.
\\

For the radial path from $u=0$ to a root $\xi_{j(d-1)}$
of $u^d=0$ the degeneration is characterized
by the following properties:
\begin{enumerate}
\item
the order of the punctures according to argument is preserved,
which follows from the provable fact, that two punctures never have the
same argument, ie. belong never to the same radial ray,
\item
the punctures with index congruent to $j \mod d$ converge to infinity,
\item
the trace of all other punctures remains bounded.
\end{enumerate}

When $u$ turns in a small circle around $\xi_{j(d-1)}$
most punctures move but very little.
In contrast
the punctures close to infinity turn on a large circle by the $(d-1)^{\text{th}}$
part of the full circle.
With the radial contraction of $u$ they retrace the movement of their 
$d^{\text{th}}$ neighbour.

Any braid thus obtained is 
$$
\s_1 \cdots \s_{(d-1)d-1}\,
\s_{(d-1)d-1} \cdots \s_{(d-2)d+1}\,
\s_{(d-2)d-1} \cdots \s_{(d-3)d+1}\,\cdots\,
\s_{d-1} \cdots \s_1,
$$
or a conjugate of it by a power of $\s_{(d-1)d-1}\cdots\s_1$.

We are free to modify these braids by an element from Prop. \ref{hurwitz},
and we use this freedom to multiply with braids of the form
$\s_{j}^{-1} \cdots \s_{j+d-2}^{-1}\s_{j+d-1}\s_{j+d-2} \cdots\s_{j}$ to
get 
$$
\s_1 \cdots \s_{(d-1)d-1}\,
\s_{(d-1)d-1} \cdots \s_1,
$$
and its conjugates by $\s_{(d-1)d-1}\cdots\s_1$.

These braids have to be conjugated as in the proof of Thm. \ref{bbmbp} to
get braids which fit with the braids of Thm. \ref{bbmbp}.
To get to our claim we have thus to modify again.
\eqed

\section{Braid monodromies versus geometric monodromies}
\label{geo}

The topological analogue of a plane curve with simple branching
along a preferred projection is a simply braided surface:

\begin{definition}
A \emph{simply braided surface} is a submanifold of dimension two
with boundary
$$
(S,\partial S) \quad \subset \quad (D^2\times \CC, S^1 \times \CC)
$$
such that
\begin{enumerate}
\item
The induced projection $S \to D^2 $ is a simple branched covering.
\item
The induced projection $\del S \to S^1$ is an unbranched covering.
\end{enumerate}
\end{definition}

\begin{remark}
By the Riemann-Hurwitz formula,
the \emph{branch set} $\bfrak$ of branch points of $(S,\del S)$ has 
cardinality $|\bfrak| = d - e( S ) $, where $d$ is the degree of the covering.
\pagebreak[4]
\end{remark}

Accordingly the range of our geometric monodromy will
consist of mapping classes preserving the braided surface structure.

\begin{definition}
The \emph{braided mapping class group} $\mfami(S)$ is the group of
isotopy classes of orientation preserving
diffeomorphisms of $D^2\times \CC$ which
\begin{enumerate}
\item
preserve $(S,\del S)$,
\item
permute the fibres of $D^2\times \CC\to D^2$,
\item
preserve the fibres of  $S^1\times \CC\to S^1$,
\item
are compactly supported.
\end{enumerate}
\end{definition}

In the case of versal unfoldings rel $\prx $ of bivariate polynomials
the zero set of a generic polynomial is naturally a simply braided
surface.

In the case of projective plane curves we look at curves not in the
degeneracy locus. The intersection of such a curve
with the complement of a small tubular neighbourhood of the line at
infinity determines a simply braided surface up to isomorphism.
The same is true for families of such curves over a loop.
The family of boundaries need not be trivial.

In any case
we may define:

\begin{definition}
The \emph{braided geometric monodromy}
is defined on the fundamental group of the
complement of the degeneracy locus and takes
values in the braided mapping class group of
\begin{enumerate}
\item
zero set of a generic polynomial in the case of bivariate polynomials,
\item
complement of a tubular neighbourhood of the line at infinity
in the case of generic projective plane curves.
\end{enumerate}
\end{definition}

Since the magnitude of $D^2$ can be chosen large in comparison with
the deformation parameters, we may deduce that the family of 
boundaries is trivialisable in the case of generic bivariate polynomials.

This observation serves well in the proof of the following comparison
result.

\begin{proposition}
\label{isom1}
The braided geometric monodromy group of versal unfolding rel $\prx $
of a bivariate polynomial with isolated singularity is isomorphic to
the bifurcation braid monodromy group.
\end{proposition}

\begin{proposition}
\label{isom}
Given a versal unfolding rel $\prx $ of a plane curve singularity.
Then the following two monodromy groups are isomorphic:
\begin{enumerate}
\item
the bifurcation braid monodromy group,
\item
the braided geometric monodromy group.
\end{enumerate}
\end{proposition}

\proof
A representative of a braided mapping class induces a diffeomorphism
of the base $D^2$ preserving the singular values $\bfrak$.
The braided mapping class thus determines a mapping class of the
punctured base.
Hence the bifurcation braid monodromy map factors through
braided geometric monodromy map. In fact the bifurcation braid
is naturally identified with the induced mapping class of
$(D^2,\bfrak)$.

Conversely we note that a braid in the bifurcation braid monodromy
determines a unique braided mapping class. On one hand
it fixes an induced mapping class on $(D^2,\bfrak)$. On the other hand
the map on the boundary is trivial.
Hence there is a unique lift to the braided mapping class group.
\eqed

\pagebreak[3]

\begin{remark}
The same result is true in the case of plane projective curves,
but the proof is more involved, since one has to take into account
that the map on the boundary may vary.
The crucial step is in fact to determine the map on the boundary
from the braid.
\end{remark}

There are obvious maps from the braided 
geometric monodromy of projection germs
to the knotted geometric monodromy of plane curve singularities and further
to the geometric monodromy.

\begin{proposition}
The knotted geometric monodromy is injective for plane curve
singularities of type $A_n$ and $D_n$.
\end{proposition}

\proof
This follows immediately from \cite{pv}, since geometric monodromy
factors through knotted geometric monodromy.
\eqed

We know of the failure of the geometric monodromy to be injective
in general by the result of Wajnryb \cite{w}. 
We also know of its failure to be surjective, see the result of Hirose
\cite{hirose} in the case of projective plane curves.

But there is hope that knotted geometric monodromy
is better in the sense that injectivity and surjectivity hold true
or fail at least to a lesser extend.

\begin{acknowledgement}
The author would like to thank the Mathematische Institut G\"ottingen, where
this work was completed.
He gratefully acknowledges the support of the Forschergruppe 790 ÔClassification of algebraic surfaces and compact complex manifoldsÕ of the DFG (Deutsche Forschungsgemeinschaft).
\end{acknowledgement}

\end{document}